\documentclass{amsart}
\usepackage{amssymb,amsmath,amsthm,xspace}
\usepackage[small]{diagrams}

\def\mat#1{\ensuremath{#1}\xspace}

\def\cF{\mat{\mathbb{F}}}
\def\cN{\mat{\mathbb{N}}}   
\def\cQ{\mat{\mathbb{Q}}}   

\def\cZ{\mat{\mathbb{Z}}}   

\def\lE{\mat{\mathcal{E}}}

\def\lH{\mat{\mathcal{H}}}

\def\lR{\mat{\mathcal{R}}}

\let\Phitemp\Phi \def\Phi{\mat{\Phitemp}}
\let\Psitemp\Psi \def\Psi{\mat{\Psitemp}}
\let\etatemp\eta \def\eta{\mat{\etatemp}}

\def\la{\mat{\lambda}}

\let\mutemp\mu \def\mu{\mat{\mutemp}}
\let\nutemp\nu \def\nu{\mat{\nutemp}}
\def\si{\mat{\sigma}}

\def\al{\mat{\alpha}}
\def\be{\mat{\beta}}
\def\ga{\mat{\gamma}}
\def\Ga{\mat{\Gamma}}
\def\de{\mat{\delta}}

\let\xitemp\xi \def\xi{\mat{\xitemp}}

\def\te{\mat{\theta}}

\def\mrm@#1{\mat{\mathrm{#1}}}

\def\gm{\mat{\mathfrak{m}}}

\def\DMO{\DeclareMathOperator}

\DMO{\Hom}{Hom}
\DMO{\lHom}{\lH\mathit{om}}
\DMO{\Ext}{Ext}
\DMO{\lExt}{\lE\mathit{xt}}
\DMO{\End}{End}
\DMO{\Aut}{Aut}
\DMO{\Fun}{Fun}
\DMO{\Tor}{Tor}
\DMO{\ext}{ext}

\DMO{\Ob}{Ob}
\DMO{\Mor}{Mor}
\DMO{\im}{im}
\DMO{\coim}{coim}
\DMO{\coker}{coker}
\DMO{\Arr}{Arr}

\DMO{\Id}{Id}
\DMO{\add}{add} 
\DMO{\ind}{ind} 
\DMO{\pro}{pro} 
\DMO{\Map}{Map} %
\DMO{\Iso}{Iso} %
\DMO{\Isom}{Isom}%
\DMO{\Ind}{Ind}

\DMO{\Presh}{Presh}

\DMO\coalg{Coalg}
\DMO{\Rep}{Rep}
\DMO{\Cor}{Cor}

\DMO{\Mod}{Mod}
\DMO{\rad}{rad}
\DMO{\soc}{soc}
\DMO{\ann}{ann}

\DMO{\Spec}{Spec}
\DMO{\spec}{Spec}
\DMO{\Proj}{Proj}
\DMO{\supp}{supp}
\DMO{\Coh}{Coh}
\DMO{\coh}{Coh}
\DMO{\Qcoh}{QCoh}
\DMO{\QCoh}{QCoh}
\DMO{\Pic}{Pic}
\DMO{\Div}{Div}
\DMO{\ch}{ch}
\DMO{\Hilb}{Hilb}
\DMO{\Fitt}{Fitt}
\DMO{\Quot}{Quot}

\DMO{\Gras}{Gr}
\DMO{\Flag}{Flag}

\DMO{\cone}{cone}
\DMO{\Tw}{Tw}
\DMO{\rank}{rk}
\DMO{\rk}{rk}
\DMO{\codim}{codim}
\DMO{\cov}{cov}
\DMO{\sgn}{sgn}
\DMO{\td}{td}
\DMO{\GL}{GL}
\DMO{\SL}{SL}

\DMO\Der{Der}
\DMO\der{Der}
\DMO\coder{Coder}
\DMO{\diag}{diag}
\DMO{\HMod}{HMod} 
\DMO{\ad}{ad}
\DMO*{\colim}{colim}
\DMO*{\hocolim}{hocolim}
\DMO*{\holim}{holim}
\DMO{\Ho}{Ho}

\DMO{\har}{char}
\DMO{\sk}{sk}
\DMO{\cosk}{cosk}
\DMO{\Gal}{Gal}
\DMO{\tr}{tr}
\DMO{\Tr}{Tr}
\DMO{\Sh}{Sh}
\DMO{\Is}{Is} 
\DMO{\Hol}{Hol} 
\DMO{\Lie}{Lie} 
\DMO{\Res}{Res} 
\DMO{\irr}{irr} %
\DMO{\Irr}{Irr} %
\DMO{\Exp}{Exp} %
\DMO{\Log}{Log} %
\DMO{\Pow}{Pow}
\DMO{\mult}{mult} %
\DMO{\height}{ht} %
\DMO{\wt}{wt}
\DMO{\Vect}{Vect}
\DMO{\moda}{mod}

\def\iso{\simeq}

\def\sb{\subset}

\def\xx{\times}

\def\ms{\backslash} 
\def\pser#1{[\![#1]\!]} 
\def\lser#1{\mbox{(\!(}#1\mbox{)\!)}} 

\def\inv{^{-1}}

\def\tw{^{\rm tw}}

\def\cond#1#2{\left.#1\right|_{#2}}  

\def\ub#1{\mat{\overline{#1}}}  

\def\ang#1{\mat{\left\langle #1\right\rangle}}

\def\set#1{\mat{\{ #1\}}}

\def\arrowsUsual{
\newarrow{TeXto}----{->}
\newarrow{TeXinto}C---{->}
\newarrow{TeXonto}----{->>}
\def\ar{\rightarrow}
\def\emb{\hookrightarrow}
\def\mto{\mapsto}
\def\arr{rTeXto}
\def\embb{\rTeXinto}
\newarrow{Eq}=====
    }

\newif\ifukr\ukrfalse
\newif\ifrus\rusfalse
\newif\ifger\gerfalse

\def\theorems{
\newtheorem{nthr}{Auxiliary}[section] 
\newtheorem{prp}[nthr]{Proposition}
\newtheorem{thr}[nthr]{Theorem}
\newtheorem{lmm}[nthr]{Lemma}
\newtheorem{crl}[nthr]{Corollary}
\newtheorem{clm}[nthr]{Claim}

\theoremstyle{definition}
\newtheorem{dfn}[nthr]{Definition}
\newtheorem{rmr}[nthr]{Remark}
\newtheorem{exm}[nthr]{Example}
}

\def\theoremsElsevier{

\newtheorem{conj    [nthr]{Conjecture}
\theoremstyle{definition}

}
\theorems\arrowsUsual

\begin{document}
\title[]{On the number of stable quiver representations over finite fields}%

\author{Sergey Mozgovoy}%
\author{Markus Reineke}%


\email{mozgov@math.uni-wuppertal.de}%
\email{reineke@math.uni-wuppertal.de}%
\thanks{}%

\keywords{}%

\begin{abstract}
We prove a new formula for the generating function of polynomials 
counting absolutely stable representations of quivers over finite fields.
The case of irreducible representations is studied in more detail.
\end{abstract}
\maketitle

\section{Introduction}
Many objects arising in representation-theoretic contexts have the favourable property that they are counted over finite fields by
evaluating polynomials with integer coefficients at the cardinality of the field. In particular, this applies to the rational points of many moduli spaces of representations, as well as other varieties appearing in representation theory (see \cite[Appendix]{HAK} or \cite{vdBE} for a general discussion of varieties whose rational points over finite fields are counted by polynomials).

This special property opens the very interesting possibility to prove, or at least predict, topological and geometric properties of such varieties from the polynomials arising in this way. For example, the Euler characteristic of such varieties can be computed as the evaluation of the counting polynomial at $1$. 
It is also typical that constant term conjectures or positivity conjectures on the polynomials can be formulated, often with very interesting conjectural geometric interpretations. A famous example is provided by the Kac conjecture, which interprets the constant terms of the polynomials counting absolutely indecomposable quiver representations of fixed dimension vector as root multiplicities in corresponding Kac-Moody algebras (this conjecture was proved in \cite{CBVdB} in the case of indivisible dimension vectors; a proof in the general case was announced in \cite{Ha}). Furthermore, the mere existence of these polynomials might substantiate conjectures on rationality
of varieties or the existence of cell decompositions.
  
\smallskip  
The second named author proved the polynomiality property for the number of isomorphism classes of absolutely stable (in particular,  irreducible) quiver representations, equivalently, rational points over finite fields of moduli of stable quiver representations in \cite{Rei1}, using Hall algebra techniques and the arithmetic of representations.
However, the resulting formula for these polynomials is of a recursive nature, and it is highly complicated. No special properties of the polynomials could be extracted, except for some predictions based on computer experiments.

The first named author already established the usefulness of $\lambda$-ring techniques for the study of counting polynomials in \cite{M2}, thereby deriving a combinatorial reformulation of the Kac conjecture mentioned above.

\smallskip
In the present paper both approaches are combined, by implementing $\lambda$-ring techniques into the methods of \cite{Rei1}. Our first result, Theorem \ref{thr:main}, provides a simple formula for the generating function of the polynomials counting absolutely stable quiver representations, in terms of certain rational functions introduced in \cite{Rei2} (and also used in \cite{ER}). Our second result, Theorem \ref{thr:main irreducible},
gives a further simplification of the formula when the 
irreducible representations are counted.

In the latter case, computer experiments suggest to consider the expansion of the counting polynomials at $q=1$. In \cite{Rei1}, an explicit formula was conjectured for the linear term of this expansion (the constant term being zero in general). This conjecture was proved in \cite{Rei3} using entirely different methods, namely, fixed point localization in quiver moduli. Our third result, Theorem \ref{thr:recursive irred} and its corollary, give a new formula for the linear term, and rederive the result of \cite{Rei3} in the case of multiple loop quivers.

Finally, we conjecture that the Taylor coefficients around $q=1$ are nonnegative. Based on further computer experiments, a conjecture on these Taylor coefficients is formulated in Remark \ref{finalremark}. We hope that a proof of this conjecture, and a combinatorial interpretation of the Taylor coefficients, could give detailed information on the geometry of moduli spaces of irreducible quiver representations.

\smallskip
The paper is organized as follows: Section \ref{sec:combi stuff} discusses the relevant $\lambda$-ring methods, and Section \ref{sec:repres stuff} collects the necessary facts on (stable) quiver representations and Hall algebras. The first main result, the formula for the generating function of polynomials counting absolutely stable representations, is derived in Section \ref{sec:computation}, and specialized to the case of irreducible representations in Section \ref{sec:irreducible}. The final Section \ref{sec:recursive} gives a recursive formula for the latter case and rederives the formula for the linear term of the counting polynomials 
in the case of multiple loop quiver.
\section{Combinatorial prerequisites}\label{sec:combi stuff}
We endow the ring of power series $R=\cQ(q)\pser{x_1,\dots,x_r}$
with a structure of a \la-ring (see \cite[Appendix]{M2}) in terms of 
Adams operations by
$$\psi_n(f(q,x_1,\dots,x_r))=f(q^n,x_1^n,\dots,x_r^n).$$
As a \la-ring, the ring $R$ possesses
also \la-operations and \si-operations.
Let \gm be the unique maximal ideal of the ring
$R$. We define the map
$\Exp:\gm\ar 1+\gm$ by the formula 
$$
\Exp(f)=\sum_{k\ge0}\si_k(f)=\exp\left(\sum_{k\ge1}\frac 1k\psi_k(f)\right).
$$
This map has an inverse $\Log:1+\gm\ar\gm$, given by
the Cadogan formula \cite{Cadogan1,Getz1,M2}
$$
\Log(f):=\sum_{k\ge1}\frac{\mu(k)}k\psi_k(\log(f)),
$$
where $\mu$ is a classical M\"obius function.
We define the map $\Pow:(1+\gm)\xx R\ar 1+\gm$ by
the formula
$$
\Pow(f,g):=\Exp(g\Log(f)).
$$
The usual power map is defined by $f^g=\exp(g\log(f))$ for
$f\in 1+\gm$ and $g\in R$.

The following useful technical result relates the map 
$\Pow$ with the usual power map.

\begin{lmm}[{\cite[Lemma 22]{M2}}]\label{lmm:power formula}
Let $f\in 1+\gm$ and $g\in R$. 
Define the elements $g_d\in R$, $d\ge1$ by the formula 
$\sum_{d\mid n}d\cdot g_d=\psi_n(g)$, $n\ge1$. Then we have
$$\Pow(f,g)=\prod_{d\ge1}\psi_d(f)^{g_d}.$$
\end{lmm}

\begin{rmr}
The maps $\Exp$, $\Log$ and $\Pow$ can be defined
for an arbitrary complete \la-ring, see \cite[Appendix]{M2}.
Their basic properties can be also found there.
\end{rmr}

For $n\in\cZ$ and $m\ge0$, define the q-binomial coefficients (see \cite{M4})
$$[n,m]=\frac{\prod_{i=1}^m(1-q^{n+i})}{\prod_{i=1}^m(1-q^{i})}\in \cQ(q),$$
$$[\infty,m]=\frac 1{\prod_{i=1}^m(1-q^i)}\in \cQ(q).$$
Let $I=\set{1,\dots,r}$ be a finite set.
For any $\la\in\cZ^I$ and $\al\in\cN^I$, define
$$[\la,\al]:=\prod_{i\in I}[\la^i,\al^i],\qquad 
[\infty,\al]:=\prod_{i\in I}[\infty,\al^i].$$
Define
$$p=p(x_1,\dots,x_r)
=\sum_{\al\in\cN^I}[\infty,\al]x^\al\in \cQ(q)\pser{x_1,\dots,x_r}.$$
For any $\la\in\cZ^I$, define
$$p(\la)=p(\la;x_1,\dots,x_r)
=\sum_{\al\in\cN^I}[\la,\al]x^\al\in\cQ(q)\pser{x_1,\dots,x_r}.$$

As a particular case of the Heine formula \cite[Theorem 13.1]{KC1}, 
one can show for $r=1$ 

\begin{lmm}[see \cite{M4}]\label{lmm:formula for p}
We have
$$p(x)=\Exp\left(\frac x{1-q}\right)$$
and, for any $n\in\cZ$,
$$p(n;x)=\Exp\left(\frac{1-q^{n+1}}{1-q}x\right).$$
\end{lmm}

This lemma implies that
$$
p(x_1,\dots,x_r)
=\prod_{i\in I}\left(\sum_{k\ge0}[\infty,k]x_i^k\right)
=\prod_{i\in I}p(x_i)
=\Exp\left(\frac{\sum_{i\in I}x_i}{1-q}\right).
$$
The functions $p$ and $p(\la)$ can be related as follows.
For any $f\in R$, define the conjugation
$\ub{f(q)}:=f(q\inv)$. 
For any $\la\in\cZ^I$,
define the algebra homomorphism $S_\la:R\ar R$ by
$S_\la(x^\al)=q^{\sum \la^i\al^i}x^\al$. 

\begin{lmm}[see \cite{M4}]\label{lmm:rel1}
For any $\la\in\cZ^I$, we have $p(\la)=p\cdot S_\la(\ub p)$.
\end{lmm} 
\section{Quiver representations}\label{sec:repres stuff}
Let $(\Ga,I)$ be a finite quiver. Here $\Ga$ denotes the set of
arrows and $I$ denotes the set of vertices of the quiver. Usually we denote
the quiver just by \Ga. 
For any field $k$, we denote by $k\Ga$
the path algebra of \Ga over $k$. Any representation $M$ of $k\Ga$ can
be decomposed into a sum of $k$-vector spaces 
$M=\oplus_{i\in I}M_i$ in a natural way. We define
the dimension of $M$ to be
$\dim M=(\dim M_i)_{i\in I}\in\cN^I$.

We define the Ringel form on $\cZ^I$ to be the  bilinear form 
defined by 
$$\ang{\al,\be}=\sum_{i\in I}\al^i\be^i-\sum_{h:i\ar j}\al^i\be^j,$$
for any $\al=(\al^i)_{i\in I}\in\cZ^I$ and 
$\be=(\be^i)_{i\in I}\in\cZ^I$. Define the Tits form on
$\cZ^I$ to be the quadratic form on $\cZ^I$ defined by
$T(\al)=\ang{\al,\al}$, for any $\al=(\al^i)_{i\in I}\in\cZ^I$.

Define the space of $k\Ga$-representations of dimension $\al\in\cN^I$
to be
$$\lR_\al(k)=\bigoplus_{h:i\ar j}\Hom_k(k^{\al^i},k^{\al^j}).$$
The group $\GL_\al(k)=\prod_{i\in I}\GL_{\al^i}(k)$ acts on 
$\lR_\al(k)$ by
$$(g_i)_{i\in I}\cdot (x_h)_{h\in\Ga}=(g_{j}x_h g_{i}\inv)_{h : i\ar j}$$
and the orbits of this action parametrize the isomorphism classes
of $k\Ga$-rep\-re\-sen\-ta\-tions of dimension \al.

\subsection{Stable representations}
We define a stability map (or just stability) to be a map $\te:I\ar\cZ$. It can be extended
by linearity to $\te:\cZ^I\ar\cZ$.
Associated with $\te$, we define a slope function 
$\mu:\cN^I\ms\set0\ar\cQ$ by
$$\mu(\al)=\frac{\te(\al)}{\height\al},$$
where $\height\al=\sum_{i\in I}\al^i$. For any nonzero
$k\Ga$-representation $M$, we define $\mu(M)=\mu(\dim M)$.

Let us fix once and for all some stability \te.
A $k\Ga$-representation $M$ is called
stable (respectively, semistable) if for any nonzero proper 
subrepresentation $N\sb M$, we have
$\mu(N)<\mu(M)$ (resp. $\mu(N)\le\mu(M)$). 
A stable representation of $k\Ga$ is called absolutely stable if it stays stable
after any  finite field extension of $k$.
A semistable representation is called polystable if it is a direct sum of stable
representations.

For any $\mu\in\cQ$, we denote
by $\moda_\mu k\Ga$ the category of semistable
$k\Ga$-rep\-re\-sen\-ta\-tions having slope $\mu$ 
(by convention, the zero representation is also contained in this category).
The category $\moda_\mu k\Ga$ is an abelian category and its simple 
objects are precisely the stable representations. 
The subcategory of $\moda_\mu k\Ga$ 
consisting of stable representations is denoted by $\moda^s_\mu k\Ga$.

Let $\al\in\cN^I$. Define $\lR^{ss}_\al(k)\sb\lR_\al(k)$ to be
the subset of semistable representations. It is invariant under
the action of $\GL_\al(k)$ and the orbits
of this action parametrize the isomorphism classes
of semistable representations of dimension \al. 
For any finite field $\cF_q$, we denote by $r_\al(q)$
the number
\begin{equation}\label{eq:definition of r}
r_\al(q)=\frac{\#\lR_\al^{ss}(\cF_q)}{\#\GL_\al(\cF_q)}.
\end{equation}

It was proved in \cite{Rei2} that the $r_\al(q)$ are rational
functions in the variable $q$. More precisely, we have

\begin{prp}[see {\cite[Corollary 6.2]{Rei2}}]
For any $\al\in\cN^I\ms\set0$, we have
$$r_\al(q)
=\sum_{(\al_1,\dots,\al_k)}
(-1)^{k-1}q^{-\sum_{i<j}\ang{\al_i,\al_j}}
\prod_{i=1}^k\frac{\#\lR_{\al_k}(\cF_q)}{\#\GL_{\al_k}(\cF_q)},$$
where the sum runs over all tuples $(\al_1,\dots,\al_k)$
of vectors in $\cN^I\ms\set0$, such that $\sum_{i=1}^k\al_i=\al$
and $\mu(\sum_{i=1}^j\al_i)>\mu(\al)$ for any $1\le j<k$.
\end{prp}

The goal of the paper is to provide a simple formula expressing
the numbers of absolutely stable representations
in terms of the rational functions $r_\al$.

\subsection{Absolutely stable representations}
Given a $k\Ga$-representation $S$, we define 
$r_S=\dim_k\End(S)$. A stable representation
$S$ is absolutely stable if and only if $r_S=1$, 
see \cite[Proposition 4.4]{Rei1}. 

Let $\al\in\cN^I$ and $r\ge1$.
Given a finite field $\cF_q$, we denote by
$a_\al(q)$ the number of isomorphism classes of absolutely stable 
$\cF_q\Ga$-representations of dimension \al.
We denote by $s_{\al,r}(q)$ the number of isomorphism classes of stable 
$\cF_q\Ga$-representations $S$ of dimension $\al$ and with $r_S=r$.
It was proved in \cite{Rei1} that $a_\al(q)$ 
and $s_{\al,r}(q)$ are polynomials in $q$, the first one having 
integer coefficients. 
These polynomials are related with each other as follows

\begin{lmm}[cf. {\cite[Proposition 4.5]{Rei1}}]
For any $\al\in\cN^I$ and any $r\ge1$, we have
$$\psi_r(a_\al)=\sum_{k\mid r}k s_{k\al,k}.$$
\end{lmm}
\begin{proof}
According to \cite[Proposition 4.5]{Rei1}, we have
$$s_{\al,r}(q)=\frac 1r\sum_{k\mid r}\mu(r/k)a_{\al/r}(q^k),$$
Using the Adams operations, we can rewrite this formula
in the form
$$rs_{r\al,r}=\sum_{k\mid r}\mu(r/k)\psi_k(a_{\al}).$$
The M\"obius inversion formula implies then the statement of
the lemma. 
\end{proof}

Applying the above lemma together with Lemma \ref{lmm:power formula}
we get

\begin{crl}\label{crl:power a}
Let $\gm$ be the maximal ideal of $\cQ(q)\pser{x_i|i\in I}$
and let $f\in 1+\gm$. Then
$$\Pow(f,a_\al)=\prod_{r\ge1}\psi_r(f)^{s_{r\al,r}}.$$
\end{crl}

\subsection{Hall algebra}
Let $k$ be a finite field. One defines the Hall algebra (see \cite{Ringel1})
$H(k\Ga)$
to be a $\cQ$-vector space with a basis consisting of isomorphism
classes of $k\Ga$-representations. For any 
two representations $M,\,N$ of $k\Ga$ the multiplication
$[M]\cdot[N]$ is defined by
$$[M]\cdot[N]=\sum_{[X]}g_{M,N}^X[X],$$
where $g_{M,N}^X$ denotes the number of subrepresentations
$U\sb X$ such that $U\iso N$ and $X/U\iso M$.

The Hall algebra $H(k\Ga)$ possesses a grading given by the height
of the dimension of representations. The corresponding completion
of $H(k\Ga)$ is denoted by $H\lser{k\Ga}$. It is a local algebra and
its invertible elements are those having nonzero component of degree
zero.

Endow the algebra $\cQ(q)[x_i|i\in I]$ with a new product defined by 
$$x^\al\circ x^\be=q^{-\ang{\al,\be}}x^{\al+\be}.$$
The new algebra is denoted by $\cQ(q)\tw[x_i|i\in I]$.
Given a finite field $\cF_q$, we endow the algebra $\cQ[x_i|i\in I]$
with a new product by the same formula as above, with $q$ considered
now as a rational number. 
The new algebra is denoted by $\cQ_{\cF_q}\tw[x_i|i\in I]$
or just by $\cQ_q\tw[x_i|i\in I]$.
The algebra $\cQ(q)\tw[x_i|i\in I]$ (resp. $\cQ_q\tw[x_i|i\in I]$)
has a grading given by $\deg x^\al=\height\al$. The corresponding
completion is denoted by 
$\cQ(q)\tw\pser{x_i|i\in I}$ (resp. $\cQ_q\tw\pser{x_i|i\in I}$). 
As a vector space it is
isomorphic to $\cQ(q)\pser{x_i|i\in I}$ (resp. $\cQ\pser{x_i|i\in I}$).

\begin{lmm}[see {\cite[Lemma 3.3]{Rei1}}]\label{lmm:integral}
Let $k$ be a finite field. 
Then the map $\int:H(k\Ga)\ar\cQ_k\tw[x_i|i\in I]$
given by
$$\int[M]=\frac1{\#\Aut M}x^{\dim M}$$
is a homomorphism of graded algebras. It induces
a homomorphism of completions 
$\int:H\lser{k\Ga}\ar\cQ_k\tw\pser{x_i|i\in I}$.
\end{lmm}

Let $\mu$ be a rational number and let $k$ be a finite field.
Define the element
$$e_\mu=e_{\mu,k}=\sum_{[M]\in\moda_\mu k\Ga}[M]$$
in the completed Hall algebra $H\lser{k\Ga}$.
Its component of degree zero 
is nonzero and therefore $e_\mu$ is invertible.

\begin{lmm}[see {\cite[Lemma 3.4]{Rei1}}]\label{lmm:e inverse}
We have
$$e_{\mu,k}\inv=\sum_{[M]\in \moda_\mu k\Ga}\ga_M[M]$$
where 
$$\ga_M=
\begin{cases}
0& 
\text{if }M\text{ is not polystable},\\
\prod_{[S]\in\moda^s_\mu k\Ga}(-1)^{m_S}(\#\End S)^{\binom{m_S}2}&
\text{if }M=\bigoplus_{[S]\in\moda^s_\mu k\Ga} S^{m_S}.
\end{cases}
$$ 
\end{lmm}

\begin{lmm}[see {\cite{Rei1}}]\label{lmm:integral of e}
Let $k$ be a finite field. Then
$$\int e_{\mu,k}=\sum_{\al\in\cN^I}r_\al(\#k)x^\al.$$
\end{lmm}

In the next section, we will express 
$\int e_\mu^{-1}$ in terms of the polynomials $a_\alpha$, drastically 
simplifying the corresponding formula from \cite{Rei1}.
This will allow us to relate the polynomials $a_\al$ with the rational
functions $r_\al$ and to get a nice formula for the computation
of $a_\al$. 
\section{Computation of the integral}\label{sec:computation}
Let $(\Ga,I)$ be a finite quiver, $\te:I\ar\cZ$ be a stability
and $\mu$ be a rational number. In this section, 
the slope function and the
stability condition are always determined by \te.
Representations are assumed to have slope \mu, if not stated
otherwise.

Let us define the generating functions 
$$a(q)=\sum_{\al\in\cN^I}a_\al(q) x^\al,\qquad r(q)=\sum_{\al\in\cN^I}r_\al(q) x^\al$$
in $\cQ(q)\pser{x_i|i\in I}$, where the polynomials $a_\al\in\cZ[q]$
and the rational functions $r_\al\in\cQ(q)$ were defined in 
Section \ref{sec:repres stuff}.

\begin{thr}\label{thr:main}
We have
$$r(q)\circ \Exp\left(\frac{a(q)}{1-q}\right)=1$$
in the algebra $\cQ(q)\tw\pser{x_i|i\in I}$.
\end{thr}
\begin{proof}
It is enough to prove the corresponding formula for the specializations
of $q$ to the numbers of points of finite fields, i.e.,
for any finite field $k$, we have to prove
$$\cond{r(q)}{q=\#k}\circ\cond{\Exp\left(\frac{a(q)}{1-q}\right)}{q=\#k}=1$$
in $\cQ_k\tw\pser{x_i|i\in I}$.

By Lemma \ref{lmm:integral of e}, we have $\int e_{\mu,k}=\cond{r(q)}{q=\#k}$.
By Lemma \ref{lmm:integral}, we have 
$$\int e_{\mu,k}\circ \int e_{\mu,k}\inv=\int(e_{\mu,k}\cdot  e_{\mu,k}\inv)=1.$$
Thus, we just have to prove
$$\int e_{\mu,k}\inv=\Exp\left(\frac{a(q)}{1-q}\right)\bigg|_{q=\#k}.$$
This is the content of the next theorem.
\end{proof}

\begin{thr}
Let $k$ be a finite field. Then
$$\int e_{\mu,k}\inv=\Exp\left(\frac{a(q)}{1-q}\right)\bigg|_{q=\#k}.$$
\end{thr}
\begin{proof}
By Lemma \ref{lmm:e inverse}, we have
$$\int e_{\mu,k}\inv=\sum_{[M]}\frac{\ga_M}{\#\Aut M}x^{\dim M}
=\sum_{(m_S)_{[S]}}\frac{\prod_{[S]}(-1)^{m_S}(\#\End S)^{\binom{m_S}2}}
{\#\Aut\left(\bigoplus_{[S]}S^{m_S}\right)}x^{\sum_{[S]}m_S\dim S}.$$
Note that
$$\Aut\left(\bigoplus_{[S]}S^{m_S}\right)\iso\prod_{[S]}\GL_{m_S}(\End S)$$
and therefore
\begin{align*}
\int e_{\mu,k}\inv
&=\sum_{(m_S)_{[S]}}\prod_{[S]}
\left(\frac{(-1)^{m_S}(\#\End S)^{\binom{m_S}2}}{\#\GL_{m_S}(\End S)}x^{m_S\dim S}\right)\\
&=\prod_{[S]}
\left(\sum_{m\ge0}\frac{(-1)^{m}(\#\End S)^{\binom{m}2}}{\#\GL_{m}(\End S)}x^{m\dim S}\right).
\end{align*}
For any stable representation $S$, the ring of endomorphisms $\End S$
is a finite field and, for any finite field $\cF_s$, we have
$$\frac{(-1)^m s^{\binom m2}}{\#\GL_m(\cF_s)}
=\prod_{i=1}^m(1-s^i)\inv=[\infty,m]\big|_{q=s}.$$
It follows
$$\int e_{\mu,k}\inv
=\prod_{[S]}\left(\sum_{m\ge0}[\infty,m]x^{m\dim S}\Big|_{q=\#\End S}\right).$$

Recall, that for any representation $S$, we have defined the integer
$r_S$ to be $\dim_k\End S$. It follows that 
$\#\End S=(\#k)^{r_S}$ and therefore
$$[\infty,m]\big|_{q=\#\End S}
=\psi_{r_S}([\infty,m])\big|_{q=\#k}.$$
This implies 
\begin{align*}
\int e_{\mu,k}\inv
&=\prod_{[S]}\left(\sum_{m\ge0}\psi_{r_S}([\infty,m])x^{m\dim S}\right)\Bigg|_{q=\#k}\\
&=\prod_{\stackrel{\al\in\cN^I}{r\ge1}}
\left(\sum_{m\ge0}\psi_{r}([\infty,m])x^{m\al}\right)^{s_{\al,r}}\Bigg|_{q=\#k}.
\end{align*}
Next
\begin{multline*}
\prod_{\al,r}\left(\sum_{m\ge0}\psi_{r}([\infty,m])x^{m\al}\right)^{s_{\al,r}}
=\prod_{\al,r}\psi_r\left(\sum_{m\ge0}[\infty,m]x^{m\al/r}\right)^{s_{\al,r}}\\
=\prod_{\al,r}\psi_r(p(x^{\al/r}))^{s_{\al,r}}
=\prod_{\al,r}\psi_r(p(x^{\al}))^{s_{r\al,r}}
\rEq^{\text{Cor. }\ref{crl:power a}}\prod_\al\Pow(p(x^\al),a_\al)\\
=\prod_\al\Exp(a_\al\Log(p(x^\al)))
\rEq^{\text{Lemma }\ref{lmm:formula for p}}\Exp\left(\sum_\al a_\al\frac{x^\al}{1-q}\right)
=\Exp\left(\frac{a(q)}{1-q}\right).
\end{multline*}

\end{proof}


\section{Irreducible representations}\label{sec:irreducible}
In this section, we take the stability map $\te:I\ar\cZ$ to be the
zero map. Then the stability of a representation is equivalent to its irreducibility
and all representations are semistable. 
The goal of this section is to give an explicit formula for the generating
function $r\in\cQ(q)\pser{x_i|i\in I}$. 
In view of Theorem \ref{thr:main}, this will give us a rather explicit
formula for the polynomials $a_\al(q)$.

Define the operator
$T$ on $\cQ(q)\pser{x_i|i\in I}$ by $T(x^\al)=q^{T(\al)}x^\al$. 
In Section \ref{sec:combi stuff} we have defined
$p=\Exp\left(\frac{\sum_{i\in I}x_i}{1-q}\right)$.
We have defined there also the conjugation $\ub {f(q)}=f(q\inv)$.

\begin{thr}\label{thr:main irreducible}
We have $r=\ub{Tp}$. In particular,
$$\ub{Tp}\circ\Exp\left(\frac{a(q)}{1-q}\right)=1$$
in the algebra $\cQ(q)\tw\pser{x_i|i\in I}$.
\end{thr}
\begin{proof}
Let $k$ be a finite field and $s=\#k$. Recall that
$$r(s)=\sum_{\al\in\cN^I}
\frac{\# \lR_\al(k)}{\#\GL_\al(k)}x^\al.$$
We have
$$\#\lR_\al(k)=s^{\al\cdot\al-T(\al)},$$
where $\al\cdot\al=\sum_{i\in I} \al^i\al^i$.
Also we have 
$$\#\GL_n(\cF_s)=\prod_{i=0}^{n-1}(s^n-s^i)=(-1)^ns^{\binom n2}\prod_{i=1}^n(1-s^i).$$
This implies 
$$\frac{s^{n^2}}{\#\GL_n(\cF_s)}
=\frac{(-1)^ns^{\binom{n+1}2}}{\prod_{i=1}^n(1-s^i)}
=\frac{1}{\prod_{i=1}^n(1-s^{-i})}
=\cond{\ub{[\infty,n]}}{q=s}.$$

We can write now
$$r(s)
=\cond{\sum_{\al\in\cN^I}q^{-T(\al)}\ub{[\infty,\al]}x^\al}{q=s}
=\cond{\ub{Tp}}{q=s}.$$
As this equality holds for all finite fields, we get
$r=\ub{Tp}$.
\end{proof}

\begin{rmr}
Assume that the quiver does not have oriented cycles. Then the irreducible
representations are precisely the one-dimensional ones
and we have $a(q)=\sum_{i\in I}x_i$. This implies
$$\Exp\left(\frac{a(q)}{1-q}\right)
=\Exp\left(\frac{\sum_{i\in I}x_i}{1-q}\right)
=p.$$
It follows that $\ub{Tp}\circ p=1$. This is a purely combinatorial
formula obtained by methods of geometric nature.
We give a combinatorial proof of it in Remark \ref{rmr:combi proof}.
\end{rmr} 
\section{Recursive formula}\label{sec:recursive}
Let $a(q)=\sum_{\al\in\cN^I}a_\al(q)x^\al$ be the generating
function for the irreducible representations of a quiver. 
In this section we will show that the expression 
$$\frac{a(q)-\sum_{i\in I}x_i}{1-q}$$
is well defined at $q=1$ and we will give a recursive formula for
the computation of both the expression and of its value at $q=1$.

Define the Ringel matrix $R$ by the formula
$\ang{\al,\be}=\al^tR\be$,
for any $\al,\,\be\in\cN^I$. Then 
$$R_{ij}=\de_{ij}-\#\set{\text{arrows from }i\text{ to }j}.$$
For any 
$f\in\cQ(q)\pser{x_i,i\in I}$, 
we denote by $f_\alpha$ the coefficient of
$x^\alpha$ in $f$.

\begin{thr}\label{thr:recursive irred}
Expression $f=\Exp\left(\frac{a(q)-\sum x_i}{1-q}\right)$ can be determined
by the recursive formula:
$$(p(-R\al)\cdot f)_\al=0$$
for any $\al>0$. This expression is well defined at $q=1$.
The specialization $f|_{q=1}$ can be determined from the same 
recursive formula, using the simplification 
$p(\la)|_{q=1}=\prod_{i\in I} (1-x_i)^{-\la^i-1}$. 
\end{thr}
\begin{proof}
We have shown that
$$\ub {Tp}\circ\Exp\left(\frac{a(q)}{1-q}\right)=1.$$
Let us denote $\Exp\left(\frac{a(q)}{1-q}\right)$ by $g$. Then,
for any $\al>0$, we have
$$\sum_{0\le\be\le\al}(\ub{Tp})_{\be}g_{\al-\be}q^{-\ang{\be,\al-\be}}=0.$$
Note that 
$$(\ub{Tp})_{\be}=q^{-T(\be)}\ub p_\be
=q^{-\ang{\be,\be}}\ub p_\be.$$
This implies
\begin{multline*}
0=\sum_{0\le\be\le\al}\ub p_\be g_{\al-\be}q^{-\ang{\be,\al}}
=\sum_{0\le\be\le\al}g_{\al-\be}q^{-\be^tR\al}\ub p_\be\\
=\sum_{0\le\be\le\al}g_{\al-\be}(S_{-R\al}\ub p)_\be
=(gS_{-R\al}\ub p)_\al.
\end{multline*}
By Lemma \ref{lmm:rel1}, we have $S_{-R\al}\ub p=p(-R\al)/p$
and therefore 
$$(p(-R\al)\cdot g/p)_\al=0$$
for any $\al>0$. 
We note that
$$g/p=\Exp\Big(\frac{a(q)-\sum x_i}{1-q}\Big)=f.$$
The above formula allows us to determine $f$ recursively.
Moreover, as the expressions $[n,k]$ are well defined at $q=1$, we 
deduce that the $p(\la)$ are well defined at $q=1$. 
The recursive formula implies that $f=g/p$ is well 
defined at $q=1$.
\end{proof}

\begin{crl}
Let $(\Ga, I)$ be a quiver with one node  and $m$ loops.
Let $a=\sum_{d\ge0}a_dx^d$ be the generating function of the polynomials
counting absolutely irreducible representation of $\Ga$ over finite fields.
Then
$$\cond{\Exp\left(\frac{a(q)-x}{1-q}\right)}{q=1}=1-mx.$$
\end{crl}
\begin{proof}
To derive the statement from our recursive formula, we have to prove that for any
$n>0$, the coefficient of $x^n$ in 
$$(1-x)^{n-mn-1}(1-mx)$$
is zero. This can be checked directly. 
\end{proof}

\begin{rmr}\label{r63}
The above corollary is equivalent to the formula
$$\frac{a_d(q)}{q-1}\Big|_{q=1}=\frac 1d\sum_{k\mid d}\mu(d/k)m^k,\qquad d\ge2.$$
proved in \cite{Rei3}. Indeed, as $a_1(q)=q^m$,
we deduce from the formula
\begin{multline*}
\frac{a(q)-q^mx}{q-1}\Big|_{q=1}+mx
=\sum_{d\ge1}\frac 1d\sum_{k\mid d}\mu(d/k)m^kx^d\\
=\sum_{k,l\ge1}\frac{\mu(l)}l\psi_l\left(\frac{m^kx^k}k\right)
=\sum_{l\ge1}\frac{\mu(l)}l\psi_l\left(\log\left(\frac1{1-mx}\right)\right)
=-\Log(1-mx).
\end{multline*}
Note that $mx=\frac{q^mx-x}{q-1}\Big|_{q=1}$. This implies
$$\frac{a(q)-x}{1-q}\Big|_{q=1}=\Log(1-mx).$$
The last formula is equivalent to the statement of the corollary.
\end{rmr}

\begin{rmr}\label{rmr:combi proof}
We can give now a combinatorial proof of the formula 
$\ub{Tp}\circ p=1$ in the case when the quiver does not have
oriented cycles. By the proof of the above theorem, this
is equivalent to the condition
$$(p(-R\al)\cdot p/p)_\al=0$$
for any $\al>0$. Thus, we have to show that
$$p(-R\al)_\al=[-R\al,\al]=0$$
for $\al>0$. We can assume that $R$ is upper-triangular
with ones on the diagonal
after permuting the nodes of the quiver. 
The above expression is a product of $q$-binomial coefficients
and it is easy to see that one of them has the form $[-n,n]$
for some $n>0$. But $[-n,n]=0$ for $n>0$.
\end{rmr}

\begin{rmr}\label{finalremark}
Expanding the polynomial $a_d(q)$ in powers of $q-1$, Remark \ref{r63} shows in particular that the linear term
$$\frac 1d\sum_{k\mid d}\mu(d/k)m^k$$
in this expansion is nonnegative (the constant term being zero except for the trivial case $d=1$): this number is well-known in combinatorics as the number of primitive necklaces consisting of $d$ beads in $m$ colours (see e.g. \cite{Reu}).

It has been verified by computer experiments, in the cases $m=2,3,4$ and $d\leq 12$, that indeed {\it all} coefficients in this expansion are nonnegative. One might therefore conjecture that
$$a_d(q)\in\cN[q-1]$$
for all $d\geq 1$.

Even more optimistically, one might hope for a geometric explanation for such a result, for example the existence of a finite stratification of the moduli space of $d$-dimensional irreducible representations (see \cite{Rei3} for the precise definition) into strata isomorphic to tori (and thus contributing powers of $q-1$ to the number of rational points of the moduli space).

The string representations of the $m$-loop quiver, introduced in \cite{Rei3} can be seen as a first step in this program: they are parametrized up to isomorphism by primitive necklaces as above, and thus provide an explanation for the linear term $\frac 1d\sum_{k\mid d}\mu(d/k)m^k$ in the above expansion.

Against this background, it is natural to study the expansion
$$\cond{\Exp\left(\frac{a(q)-x}{1-q}\right)}{q=1}=\sum_{n=0}^\infty f_n(m,t) (q-1)^n.$$
Besides our result saying that $f_0(m,t)=1-mt$, computer experiments suggest 
$$f_1(m,t)=\binom m2\frac{t(t-1)}{(1-mt)^2}$$ 
and more generally
$$f_n(m,t)\cdot(1-mt)^{3n-1}$$ 
is a polynomial in $m$ and $t$ having degree $3n-1$ with respect to $t$.
\end{rmr} 

\bibliography{a9}

\providecommand{\bysame}{\leavevmode\hbox to3em{\hrulefill}\thinspace}
\providecommand{\href}[2]{#2}
\begin{thebibliography}{10}

\bibitem{Cadogan1}
Charles Cadogan, \emph{The {M}\"obius function and connected graphs}, J. Comb.
  Theory, Ser. B \textbf{11} (1971), 193--200.

\bibitem{CBVdB}
William Crawley-Boevey and Michel Van~den Bergh, \emph{Absolutely
  indecomposable representations and {K}ac-{M}oody {L}ie algebras}, Invent.
  Math. \textbf{155} (2004), no.~3, 537--559, with an appendix by Hiraku
  Nakajima.

\bibitem{ER}
Johannes Engel and Markus Reineke, \emph{Smooth models of quiver moduli},
  \mbox{arXiv:0706.4306}.

\bibitem{Getz1}
Ezra Getzler, \emph{{Mixed Hodge structures of configuration spaces}}, Preprint
  96-61, Max-Planck-Institut f. Mathematik, Bonn,
  \mbox{arXiv:alg-geom/9510018}.

\bibitem{Ha}
Tam{\'a}s Hausel, \emph{Betti numbers of holomorphic symplectic quotients via
  arithmetic {F}ourier transform}, Proc. Natl. Acad. Sci. USA \textbf{103}
  (2006), no.~16, 6120--6124 (electronic).

\bibitem{HAK}
Tamas Hausel and Fernando Rodriguez-Villegas, \emph{Mixed hodge polynomials of
  character varieties}, \mbox{arXiv:math/0612668}, with an appendix by Nicholas
  M. Katz.

\bibitem{KC1}
Victor Kac and Pokman Cheung, \emph{Quantum calculus}, Universitext,
  Springer-Verlag, New York, 2002.

\bibitem{M2}
Sergey Mozgovoy, \emph{{A computational criterion for the Kac conjecture}},
  \mbox{arXiv:math.RT/0608321}.

\bibitem{M4}
\bysame, \emph{Fermionic forms and quiver varieties},
  \mbox{arXiv:math.QA/0610084}.

\bibitem{Rei3}
Markus Reineke, \emph{Localization in quiver moduli},
  \mbox{arXiv:math.AG/0509361}, Preprint 2005. Submitted.

\bibitem{Rei2}
\bysame, \emph{The {H}arder-{N}arasimhan system in quantum groups and
  cohomology of quiver moduli}, Invent. Math. \textbf{152} (2003), no.~2,
  349--368.

\bibitem{Rei1}
\bysame, \emph{Counting rational points of quiver moduli}, Int. Math. Res. Not.
  (2006), Article ID 70456, 19 pages.

\bibitem{Reu}
Christophe Reutenauer, \emph{Free {L}ie algebras}, London Mathematical Society
  Monographs. New Series, vol.~7, Clarendon Press, Oxford, 1993.

\bibitem{Ringel1}
Claus~Michael Ringel, \emph{Hall algebras}, Topics in algebra, Part 1 (Warsaw,
  1988), Banach Center Publ., vol.~26, PWN, Warsaw, 1990, pp.~433--447.

\bibitem{vdBE}
Theo van~den Bogaart and Bas Edixhoven, \emph{Algebraic stacks whose number of
  points over finite fields is a polynomial}, Number fields and function
  fields---two parallel worlds, Progr. Math., vol. 239, Birkh\"auser Boston,
  Boston, MA, 2005, pp.~39--49.

\end{thebibliography}
\bibliographystyle{hamsplain}

\end{document}